\documentclass[a4paper,USenglisch,cleveref,nameinlink,thm-restate]{lipics-v2021}

\usepackage{graphicx} 
\usepackage{multirow}
\usepackage{xcolor}
\usepackage{amssymb}
\usepackage{amsmath}
\usepackage{amsthm}
\usepackage{xspace}
\usepackage{graphicx}
\usepackage{complexity}
\usepackage{hyperref}
\usepackage{todonotes}
\usepackage{caption}
\usepackage{subcaption}
\usepackage[framemethod=tikz]{mdframed}
\usepackage{bibentry}
\usepackage{dsfont}
\hideLIPIcs
\nolinenumbers

\renewcommand{\R}{\ensuremath{\mathds{R}}}
\DeclareMathOperator{\ch}{ch}

\graphicspath{{figures/}}
\newcommand{\ourref}[1]{{\color{red}$\langle$\ref{#1}$\rangle$}}
\newcommand{\litcite}[1]{{\color{blue}\cite{#1}}}

\title{On $k$-planar Graphs without Short Cycles}
\author{Michael A. Bekos}{Department of Mathematics, University of Ioannina, Greece \and\url{https://myweb.uoi.gr/bekos}}{bekos@uoi.gr}{https://orcid.org/0000-0002-3414-7444}{}
\author{Prosenjit Bose}{School of Computer Science, Carleton University, Ottawa, Canada}{jit@scs.carleton.ca}{https://orcid.org/0000-0002-8906-0573}{}
\author{Aaron B{\"u}ngener}{Department of Computer Science, University of T\"ubingen, Germany}{aaron.buengener@student.uni-tuebingen.de}{}{}
\author{Vida Dujmovi{\'c}}{School of Electrical Engineering and Computer Science, University of Ottawa, Canada}{vida.dujmovic@uottawa.ca}{https://orcid.org/0000-0001-7250-0600}{}
\author{Michael Hoffmann}{Department of Computer Science, ETH Z\"urich, Switzerland \and\url{https://people.inf.ethz.ch/hoffmann/}}{hoffmann@inf.ethz.ch}{https://orcid.org/0000-0001-5307-7106}{}
\author{Michael Kaufmann}{Department of Computer Science, University of T\"ubingen, Germany \and\url{https://uni-tuebingen.de/fakultaeten/mathematisch-naturwissenschaftliche-fakultaet/fachbereiche/informatik/lehrstuehle/algorithmik/team/prof-dr-michael-kaufmann/}}{michael.kaufmann@uni-tuebingen.de}{https://orcid.org/0000-0001-9186-3538}{}
\author{Pat Morin}{School of Computer Science, Carleton University, Ottawa, Canada}{morin@scs.carleton.ca}{https://orcid.org/0000-0003-0471-4118}{}
\author{Saeed Odak}{Department of Computer Science and Electrical Engineering, University of Ottawa, Canada}{saeedodak@gmail.com}{}{}
\author{Alexandra Weinberger}{Faculty of Mathematics and Computer Science, FernUniversit\"at in Hagen, Germany
\and\url{https://www.fernuni-hagen.de/ti/en/team/alexandra.weinberger.shtml}}{alexandra.weinberger@gernuni-hagen.de}{https://orcid.org/0000-0001-8553-6661}{Part of this work was done while A.W. was at the University of Technology in Graz and supported by the Austrian Science Fund (FWF): W1230.}
\authorrunning{Bekos, Bose, B{\"u}ngener, Dujmovi{\'c}, Hoffmann, Kaufmann, Morin, Odak, Weinberger}
\Copyright{Michael~A.~Bekos, Pat~K.~Bose, Aaron~B{\"u}ngener, Vida~Dujmovi{\'c}, Michael~Hoffmann, Michael~Kaufmann, Pat~Morin, Saeed~Odak, and Alexandra~Weinberger}

\acknowledgements{This work  was  started at the Summer Workshop on Graph Drawing (SWGD~2023) in Caldana, Italy.} 

\ccsdesc[500]{Mathematics of computing~Combinatorics}
\ccsdesc[500]{Mathematics of computing~Graph theory}
\ccsdesc[500]{Human-centered computing~Graph drawings}
\keywords{Beyond-planar Graphs, k-planar Graphs, Local Crossing Number, Crossing Number, Discharging Method, Crossing Lemma}

\begin{document}

\maketitle

\begin{abstract}
We study the impact of forbidding short cycles to the edge density of~$k$-planar graphs; a \emph{$k$-planar} graph is one that can be drawn in the plane with at most~$k$ crossings per edge. Specifically, we consider three settings, according to which the forbidden substructures are $3$-cycles, $4$-cycles or both of them (i.e., girth~$\ge 5$). For all three settings and all~$k\in\{1,2,3\}$, we present lower and upper bounds on the maximum number of edges in any~$k$-planar graph on~$n$ vertices. Our bounds are of the form~$c\,n$, for some explicit constant~$c$ that depends on~$k$ and on the setting. For general~$k \geq 4$ our bounds are of the form~$c\sqrt{k}n$, for some explicit constant~$c$. These results are obtained by leveraging different techniques, such as the discharging method, the recently introduced density formula for non-planar graphs, and new upper bounds for the crossing number of~$2$-- and $3$-planar graphs in combination with corresponding lower bounds based on the Crossing Lemma.
\end{abstract}

\section{Introduction}
\label{sec:introduction}

``What is the minimum and maximum number of edges?'' is one of the most fundamental questions one can ask about a finite family of graphs. In some cases the question is easy to answer; for instance, for the class of all graphs on~$n$ vertices the answer is even trivial. Another such family is the one of planar graphs. More precisely, for planar graphs on~$n$ vertices we know from Euler's Formula that they have at most~$3n-6$ edges. Furthermore, every planar graph on~$n$ vertices can be augmented (by adding edges) to a maximal planar graph with exactly~$3n-6$ edges. Important advances have recently been made  for non-planar graphs in the context of graph drawing beyond-planarity~\cite{DBLP:books/sp/20/HT2020}. But often an answer is much harder to come by. Specifically, there exist graph classes that are relevant in Graph Drawing where exact bounds on their edge density are difficult to 
derive. 

The family of graphs that can be embedded on the Euclidean plane with at most $k$ crossings per edge, called \emph{$k$-planar}, is a notable example. Tight bounds on the edge density of these graphs, for small values of $k$, are crucial as they lead to improvements on the well-known Crossing Lemma~\cite{az-pb-10}. This was first observed by Pach and T{\'{o}}th~\cite{PachT97}, who back in 1997 presented one of the early improvements of the Crossing Lemma by introducing tight bounds on the edge density of $1$- and $2$-planar graphs. 
Since then, only two improvements emerged; one by Pach, Radoi{\v{c}}i{\'{c}}, Tardos, and T{\'{o}}th~\cite{PachRTT06} in 2006 and one by Ackerman~\cite{DBLP:journals/comgeo/Ackerman19} in 2019, both by introducing corresponding bounds on the edge density of $3$- and $4$-planar graphs, respectively. 
On the other hand, it is worth noting that these progressive refinements on the Crossing Lemma led to corresponding improvements also on the upper bound on the edge density of general $k$-planar graphs with the best one being currently $3.81\sqrt{k} n$ due to Ackerman~\cite{DBLP:journals/comgeo/Ackerman19}.
To the best of our knowledge, for $5$-planar graphs a tight bound is missing from the literature, even though it would yield 
further improvements both on the Crossing Lemma and on the upper bound of the edge density of general $k$-planar graphs.

In this work, we continue the study of this line of research focusing on special classes of graphs; in particular, on graphs not containing some fixed, so-called \emph{forbidden substructures}.
We consider three settings, according to which the forbidden substructures are $3$-cycles ($C_3$-free), $4$-cycles ($C_4$-free) or both of them (girth~$\ge 5$). For each of these settings, the problem of finding edge density bounds has been studied both in general and assuming planarity. In particular, while $C_3$-free $n$-vertex graphs may have $\Theta(n^2)$ edges, $C_4$-free graphs and graphs of girth~$5$ have at most $O(n^{\frac32})$ edges; see e.g.~\cite{DBLP:journals/jct/Furedi96,DBLP:journals/jct/Wenger91}. 
For $C_3$-free planar graphs and planar graphs of girth~$5$, one can easily derive upper bounds on their edge density using Euler's Formula; see, e.g., \cref{tab:bounds}. 
For $C_4$-free planar graphs, Dowden~\cite{DBLP:journals/jgt/Dowden16} proved that every such graph has at most $\frac{15}{7}(n-2)$ edges, and that this bound is best possible.  
For $k$-planar graphs, Pach, Spencer and T{\'{o}}th~\cite{DBLP:journals/dcg/PachST00} provided a lower bound on the crossing number of $C_4$-free $k$-planar graphs, which can be used to obtain an asymptotic upper bound of $O(\sqrt[3]{k}n)$ on the edge density of such graphs with $n$ vertices. 
Another related research branch focuses on bipartite graphs (that avoid all odd-length cycles). For this setting, Angelini, Bekos, Kaufmann, Pfister, and Ueckerdt~\cite{DBLP:conf/isaac/AngeliniB0PU18} have proposed lower and upper bounds on the edge density of several classes of graphs beyond-planarity, including $1$- and $2$-planar graphs. 

\medskip\noindent\textbf{Our contribution.} We study the class of $k$-planar graphs in the absence of $3$-cycles, $4$-cycles and both of them. Our results are summarized as follows: 

\begin{itemize}
\item For each of the aforementioned settings, we present lower and upper bounds on the maximum number of edges of~$k$-planar graphs with~$n$ vertices when~$k\in\{1,2,3\}$. Our findings are summarized in \cref{tab:bounds}. 
\item We next use these bounds to derive corresponding lower bounds on the crossing numbers of the graphs that avoid the forbidden patterns studied. %
For a summary refer to \cref{tab:crossing-lemma}. 
\item We use the two-way dependency between edge density and Crossing Lemma to derive new bounds on the edge density of $k$-planar graphs for values of $k$ greater than $3$. 
\end{itemize}

\noindent To obtain the above results, we leverage different techniques from the literature, such as the discharging method, the recently introduced density formula for non-planar graphs~\cite{DBLP:journals/corr/abs-2311-06193}, and new upper bounds for the crossing number of~$2$-- and $3$-planar graphs (\cref{thm:cr-2planar,thm:cr-3planar}) in combination with corresponding lower bounds based on the Crossing Lemma.

\captionsetup[table]{position=bottom}
\begin{table}[t]
    \centering  
    \setlength{\tabcolsep}{2pt}
    \newcommand{\spc}{\hspace{1pt}}
    \begin{tabular}{l|r|r|r|r|r|r|r|r}
       & \multicolumn{2}{c|}{unrestricted} & \multicolumn{2}{c|}{$C_3$-free} & \multicolumn{2}{c|}{$C_4$-free} & \multicolumn{2}{c}{Girth~$5$}\\    
       $k$ & lower & upper & lower & upper & lower & upper & 
       lower & upper\\\hline
       $0$   & $3n$ & $3n$ & $2n$ & $2n$ & $\frac{15n}{7}$\spc\litcite{DBLP:journals/jgt/Dowden16} & $\frac{15n}{7}$\spc\litcite{DBLP:journals/jgt/Dowden16} & $\frac{5n}{3}$ & $\frac{5n}{3}$ \\
       $1$ & $4n$\spc\litcite{bsw-bsgr-83} & $4n$\spc\litcite{bsw-bsgr-83} & $3n$\spc\litcite{DBLP:journals/dmgt/CzapPS16} & $3n$\spc\ourref{thm:1-planar-C3-upper} & $2.4n$\spc\ourref{thm:1planar-C4-lower} & $2.5n$\spc\ourref{thm:1-planar-C4-upper} & $\frac{13n}{6}$\spc\ourref{thm:1planar-girth5-lower}  & $2.4n$\spc\ourref{thm:1-planar-girth-5}\\
       $2$ & $5n$\spc\litcite{PachT97} & $5n$\spc\litcite{PachT97} & $3.5n$\spc\litcite{DBLP:conf/isaac/AngeliniB0PU18} & $4n$\spc\ourref{thm:2-planar-C3-upper} & $2.5n$\spc\ourref{thm:2planar-c4-lower} & $3.93n$\spc\ourref{thm:2planar-c4} & $\frac{16n}{7}$\ourref{thm:2planar-g5-lower} & $3.597n$\spc\ourref{thm:2planar-g5}\\
       $3$ & $5.5n$\spc\litcite{PachRTT06} & $5.5n$\spc\litcite{PachRTT06} & $4n$\spc\litcite{DBLP:conf/isaac/AngeliniB0PU18}  & $5.12n$\spc\ourref{thm:3planar-c3}  & - & 
       $4.933n$\spc\ourref{thm:3planar-c4-free}
        & $2.5n$\spc\ourref{thm:3planar-g5-lower} & $4.516n$\spc\ourref{thm:3planar-g5} \\
       \multirow{2}{*}{$k$} & \multirow{2}{*}{$\Omega(\sqrt{k})n$\spc\litcite{PachT97}}  & \multirow{2}{*}{$3.81\sqrt{k} n$\spc\litcite{DBLP:journals/comgeo/Ackerman19}} &  & $3.19\sqrt{k}n$\spc\ourref{cor:c3-k-planar-1} & &  $3.016\sqrt{k}n$\spc\ourref{cor:c4-free-k-planar-1} &  & $2.642\sqrt{k}n$\spc\ourref{cor:g5-k-planar-1}\\
       &&&&&& $O(\sqrt[3]{k})n$\spc\litcite{pach1999new} && $O(\sqrt[3]{k})n$\spc\litcite{pach1999new}\\
    \end{tabular}  
    \caption{Maximum number of edges in $k$-planar graph classes, ignoring additive constants; results from the literature are shown in blue square brackets, results form this paper are shown in red angle brackets, bounds without a citation are derived from Euler's formula.
    The lower bound for 2-planar $C_4$-free graphs trivially holds for 3-planar $C_4$-free graphs. 
    \label{tab:bounds}}
\end{table}

\captionsetup[table]{position=bottom}
\begin{table}[t]
    \centering  
    \setlength{\tabcolsep}{6pt}
    \newcommand{\spc}{\hspace{1pt}}
    \begin{tabular}{l|c|c|c|c|c}
       & \multicolumn{2}{c|}{unrestricted} & \multicolumn{1}{c|}{$C_3$-free} & \multicolumn{1}{c|}{$C_4$-free} & \multicolumn{1}{c}{Girth~$5$}\\
       Graph class & lower & upper & lower  & lower  & lower \\\hline
       $2$-planar   &  & $\frac{10n}{3}$~\ourref{thm:cr-2planar} &  & &   \\
       $3$-planar   &  & $\frac{33n}{5}$~\ourref{thm:cr-3planar}  & & &   \\
       general & $0.034\frac{m^3}{n^2}$\spc\litcite{DBLP:journals/comgeo/Ackerman19} &  & $0.049\frac{m^3}{n^2}$\spc\ourref{thm:3planar-c3} &   $0.054\frac{m^3}{n^2}$\spc\ourref{thm:2planar-c4} &   $0.071 \frac{m^3}{n^2}$\spc\ourref{thm:2planar-g5}  \\
    \end{tabular}  
    \caption{Bounds on the crossing numbers, ignoring additive constants; hold for sufficiently~large~$m$.\label{tab:crossing-lemma}}
\end{table}

\section{Preliminary Techniques and Tools}\label{sec:preliminaries}

In this section, we describe techniques that we use in our proofs, namely, the discharging method~\cite{DBLP:journals/comgeo/Ackerman19,DBLP:journals/jct/AckermanT07} (\cref{ssec:discharging}) and a method derived from a well-known probabilistic proof~\cite{az-pb-10} of the Crossing Lemma (\cref{ssec:crossinglemma}), which we formalise in the following. This section is concluded with two theorems of independent interest providing upper bounds on the number of crossings of (general) $2$- and $3$-planar graphs (\cref{ssec:crupper}).

\subsection{The Discharging Method}\label{ssec:discharging}
In some of our proofs, we employ the discharging method~\cite{DBLP:journals/comgeo/Ackerman19,DBLP:journals/jct/AckermanT07}, which is summarised as follows. Consider a biconnected graph~$G=(V,E)$ on~$|V|=n$ vertices drawn in~$\R^2$ and its planarization~$G'=(V',E')$, where at every crossing both edges are subdivided using a new vertex of degree four.
We denote the set of faces of~$G'$ by~$F'$ and call them \emph{cells}. For a face~$f\in F'$ we denote by~$\mathcal{V}(f)$ and~$\mathcal{V'}(f)$ the set of vertices from~$V$ and~$V'$, respectively, that appear on the boundary~$\partial f$ of~$f$. Furthermore, let~$|f|=|\mathcal{V'}(f)|$ denote the \emph{size} of~$f$.

To each face~$f\in F'$ we assign a charge~$\ch(f)=|\mathcal{V}(f)|+|f|-4$. Using Euler's formula~$|V'|-|E'|+|F'|=2$, it is not difficult to check (see~\cite{DBLP:journals/comgeo/Ackerman19}) that $\sum_{f\in F'}\ch(f)=4n-8$.

We then distribute these charges so as to collect a discharge of at least~$\alpha$, for some~$\alpha>0$, for every pair~$(v,f)\in V\times F'$ such that~$v\in\mathcal{V}(f)$. 

Then $4n-8=\sum_{f\in F'}\ch(f)\ge\sum_{v\in V}\alpha\deg_G(v)=2\alpha|E|$ which implies
\begin{equation}\label{eq:edge-bound}
    m=|E|\le \frac{2}{\alpha}(n-2)\,.
\end{equation}

The main challenge when applying this discharging method is to manage the redistribution of charges so that every vertex receives its due, for~$\alpha$ as large as possible. As a natural first attempt, we may have each~$f\in F'$ discharge~$\alpha$ to each~$v\in\mathcal{V}(f)$. This leaves~$f$ with a \emph{remaining} charge of
\begin{equation}\label{eq:remaining-charge}
\ch^-(f)=\ch(f)-\alpha|\mathcal{V}(f)|=(1-\alpha)|\mathcal{V}(f)|+|f|-4\,.
\end{equation}
If~$\ch^-(f)\ge 0$, for all~$f\in F'$, then we are done. However, in general, we may have~$\ch^-(f)<0$, for some~$f\in F'$. In such a case we have to find some other face(s) that have a surplus of remaining charge they can send to~$f$.

\subsection{The Crossing Lemma}\label{ssec:crossinglemma}

We can obtain upper bounds on the density also using the Crossing Lemma~\cite{acns-cfs-82}.
As a basis, we need both an upper and a lower bound for the crossing number in terms of the number of vertices and edges. Upper bounds are discussed in \cref{ssec:crupper}. In this section we derive a lower bound using the Crossing Lemma, along the lines of its well-known probabilistic proof~\cite[Chapter~40]{az-pb-10}.

\begin{theorem}\label{thm:crossing-lemma} 
  Let~$\mathcal{X}$ be a hereditary\footnote{Closed under taking induced subgraphs.} graph family and~$a,b\in\R$ such that for every~$H\in\mathcal{X}$ with~$\nu$ vertices and~$\mu$ edges we have~$\mathrm{cr}(H)\ge a\mu-b\nu$.
  Then for every graph~$G\in\mathcal{X}$ with $n$ vertices and $m$ edges with~$2am\ge 3bn$ we have 
  \[
  \mathrm{cr}(G)\ge \frac{4 a^3}{27 b^2}\cdot\frac{m^3}{n^2}\,.
  \]
\end{theorem}
\begin{proof}
  Let~$\Gamma$ be a minimum-crossing drawing of~$G$. We take a random induced subgraph~$G_p=(V_p,E_p)$ of~$G$ by selecting every vertex independently at random with probability~$p$ and consider the drawing~$\Gamma_p$ of~$G_p$ defined by~$\Gamma$. Then any such graph~$G_p$ is in~$\mathcal{X}$, and so the lower bound on~$\mathrm{cr}(G_p)$ from above holds for~$G_p$ and thus also in expectation:
  \[
  \E(\mathrm{cr}(\Gamma_p))\ge a \cdot \E(E_p)- b \cdot \E(V_p)\,.
  \]
  We have~$\E(V_p)=pn$ and~$\E(E_p)=p^2m$. Furthermore, note that~$\Gamma$ is a minimum-crossing drawing of~$G$ and, therefore, no pair of adjacent edges crosses. Thus, for a crossing to be present in~$\Gamma_p$, all four endpoints of the crossing edge pair need to be selected. Therefore, we have~$\E(\mathrm{cr}(\Gamma_p))=p^4\mathrm{cr}(\Gamma)=p^4\mathrm{cr}(G)$. Putting everything together yields
  \begin{equation}\label{eq:crossing-lemma-1}
  \mathrm{cr}(G) \ge \frac{am}{p^2} - \frac{bn}{p^3}\,.
  \end{equation}
  The function on the right hand side of the above inequality has its unique maximum at~$p=\frac{3bn}{2am}$. Setting~$p=\frac{3bn}{2am}$ to \eqref{eq:crossing-lemma-1} yields:
  \begin{equation}\label{eq:crossing-lemma-2}
  \mathrm{cr}(G)\ge \frac{4 a^3}{27 b^2}\cdot\frac{ m^3}{n^2}\,.
  \end{equation}
  As a sanity check, we need~$p\le 1$. So the bound holds for~$2am\ge 3bn$. 
\end{proof}

The simple observation that one can remove relatively few edges from a~$k$-planar graph to obtain a~$(k-1)$-planar graph allows to lift density bounds for~$i$-planar graphs, with~$i<k$, to  bounds for~$k$-planar graphs. 
By iteratively removing edges from the graph and a drawing of it with maximum number of crossings, we can show the following.

\begin{restatable}{theorem}{naivelowerbound}\label{thm:naivelowerbound}
  Let~$\mathcal{X}$ be a monotone\footnote{Closed under taking subgraphs and disjoint unions.} graph family, let~$k$ be a positive integer, and let~$\mu_i(n)$ be an upper bound on the number of edges for every $i$-planar graph from~$\mathcal{X}$ on~$n$ vertices, for~$0\le i\le k-1$. Then for every~$G\in\mathcal{X}$ with~$n\ge 4$ vertices and~$m$ edges we have
  \[
  \mathrm{cr}(G)\ge km-\sum_{i=0}^{k-1}\mu_i(n)\,.
  \]
\end{restatable}

\begin{proof}
  Consider a graph~$G\in\mathcal{X}$ and a drawing~$\Gamma$ of~$G$ with~$\mathrm{cr}(\Gamma)=\mathrm{cr}(G)$. Then, we iteratively remove an edge from the graph and the drawing that has a maximum number of crossings. As long as the number of edges in the graph is strictly greater than~$\mu_i(n)$, such an edge has at least~$i+1$ crossings. We stop when the graph is plane, with at most~$\mu_0(n)$ edges remaining. The number of crossings removed is at least
  \[
    k(m-\mu_{k-1}(n))+\sum_{i=1}^{k-1}i(\mu_i(n)-\mu_{i-1}(n))=
    km-\sum_{i=0}^{k-1}\mu_i(n)\,.\qedhere
  \]
\end{proof}

\subsection{Upper Bounds on the Crossing Number of 2- and 3-planar graphs}\label{ssec:crupper}

The Crossing Lemma provides us with pretty good lower bounds for crossing numbers. As a complement, we also need corresponding upper bounds. For a~$k$-planar graph~$G$, we have a trivial bound of~$\mathrm{cr}(G)\le km/2$. So if~$G$ is~$2$-planar, then~$\mathrm{cr}(G)\le m\le 5n-10$. But we can do better, as the following theorem demonstrates.

\begin{theorem}\label{thm:cr-2planar}
  Every $2$-planar graph on~$n\ge 2$ vertices can be drawn with at most~$(10n-20)/3$ crossings.
\end{theorem}

\begin{proof}
  Let~$G=(V,E)$ be a $2$-planar graph on~$n$ vertices, and let~$\Gamma$ be any $2$-plane drawing of~$G$ with a minimum number of crossings (among all $2$-plane drawings of~$G$). We allow multiple edges between the same pair of vertices in~$\Gamma$, but no loops nor homotopic edge pairs (that is, for each pair~$e_1,e_2$ of edges between the same two vertices, neither of the two parts of the plane bounded by the simple closed curve~$e_1\cup e_2$ is empty). Without loss of generality we assume that~$\Gamma$ is maximal $2$-plane, that is, adding any edge to~$\Gamma$ results in a graph that is not $2$-plane anymore. We may assume that adjacent edges do not cross in~$\Gamma$~\cite[Lemma~1.1]{PachRTT06}. We claim that a~$1/3$-fraction of the edges in~$\Gamma$ is uncrossed. 

  Let us first argue how the claim implies the statement of the theorem. Denote by~$x$ the number of edges that have at least one crossing in~$\Gamma$. The number~$\gamma$ of crossings in~$\Gamma$ is upper bounded by~$2x/2=x$ because every edge has at most two crossings and every crossing is formed by exactly two edges. Every $2$-planar graph on~$n\ge 3$ vertices has at most~$5n-10$ edges~\cite{PachT97}, and this bound also holds for $2$-plane multigraphs without loops or parallel homotopic edges~\cite{DBLP:conf/compgeom/Bekos0R17}. It follows that~$\gamma\le x\le\frac23(5n-10)=(10n-20)/3$.

  So it remains to prove the claim. Consider a vertex~$v$ and denote by~$X(v)$ the set of edges incident to~$v$ that have at least one crossing in~$\Gamma$. Let~$e\in X(v)$, let~$c$ denote the crossing of~$e$ closest to~$v$, let~$e^-$ denote the part of~$e$ between~$v$ and~$c$, and let~$\chi(e)$ denote the edge that crosses~$e$ at~$c$. As~$\chi(e)$ has at most two crossings, at least one of the two curves that form~$\chi(e)\setminus c$ is uncrossed. Pick such a curve and denote it by~$\chi(e)^-$. The curve~$\chi(e)^-$ has two endpoints, one of which is~$c$ and the other is a vertex of~$G$, which we denote by~$\psi(e)$. As adjacent edges do not cross in~$\Gamma$, we have~$\psi(e)\ne v$. By closely following~$e^-$ and~$\chi(e)^-$ we can draw a curve between~$v$ and~$\psi(e)$ in~$\Gamma$ that does not cross any edge of~$\Gamma$. Thus, by the maximality of~$\Gamma$ we conclude that~$v\psi(e)$ is an edge in~$\Gamma$, and it is uncrossed because~$\Gamma$ is crossing-minimal by assumption. In this way, we find an uncrossed edge~$\eta(e)=v\psi(e)$ for each~$e\in X(v)$. Different edges~$e\ne f$ in~$X(v)$ may yield the same edge~$\eta(e)=\eta(f)$. But in this case by construction~$\eta(e)=\eta(f)$ is homotopic to both~$e^-\cup\chi(e)^-$ and~$f^-\cup\chi(f)^-$, that is, the simple closed curve~$e^-\cup\chi(e)^-\cup f^-\cup\chi(f)^-$ bounds a face in~$\Gamma\setminus\eta(e)$. It follows that there is no other edge~$g\in X(v)\setminus\{e,f\}$ for which~$\eta(g)=\eta(e)$, that is, for every uncrossed edge~$u$ incident to~$v$ in~$\Gamma$ we have~$|\eta^{-1}(u)\cap X(v)|\le 2$. Therefore, at least a~$1/3$-fraction of the edges incident to~$v$ in~$\Gamma$ is uncrossed. As this holds for every vertex~$v$, it also holds globally, which completes the proof of the claim and of the theorem.
\end{proof}

%
In a similar fashion, we can obtain an improved upper bound for $3$-planar graphs, as the following theorem demonstrates. We remark that the argument used in the proof of \cref{thm:cr-3planar} does not work for larger $k > 3$.

\begin{restatable}{theorem}{crthreeplanar}\label{thm:cr-3planar}
  Every $3$-planar graph on~$n\ge 2$ vertices can be drawn with at most~$(33n-66)/5$ crossings.
\end{restatable}

\begin{proof}
 Let~$G=(V,E)$ be a $3$-planar graph on~$n$ vertices, and let~$\Gamma$ be any $3$-plane drawing of~$G$ with a minimum number of crossings (among all $3$-plane drawings of~$G$). We allow multiple edges between the same pair of vertices in~$\Gamma$, but no loops nor homotopic edge pairs. Without loss of generality we assume that~$\Gamma$ is maximal $3$-plane, that is, adding any edge to~$\Gamma$ results in a graph that is not $3$-plane anymore. We may assume that adjacent edges do not cross in~$\Gamma$~\cite[Lemma~1.1]{PachRTT06}. Denote the number of crossings in~$\Gamma$ by~$c=\mathrm{cr}(G)$. We claim that at least~$c/6$ edges are uncrossed in~$\Gamma$.

  Let us first argue how the claim implies the statement of the theorem. Denote by~$x$ the number of edges that have at least one crossing in~$\Gamma$. The number~$c$ of crossings in~$\Gamma$ is upper bounded by~$3x/2$ because every edge has at most three crossings and every crossing is formed by exactly two edges. Every $3$-planar graph on~$n\ge 3$ vertices has at most~$5.5n-11$ edges~\cite{PachRTT06}, and this bound also holds for $3$-plane multigraphs without loops or parallel homotopic edges~\cite{DBLP:conf/compgeom/Bekos0R17}. By the claim we have~$x\le 5.5n-11-c/6$ and therefore~$2c\le 3x\le 3(5.5n-11-c/6)$. It follows that~$5c/2\le 33(n/2-2)$. Rearranging terms completes the proof.

  It remains to prove the claim. We split every crossing into two so-called halfcrossings as follows. Let~$X\subset\R^2$ denote the set of all crossings of~$\Gamma$. Consider a crossing~$\chi$ of an edge~$e$ with some other edge. The \emph{halfcrossing} of~$e$ at~$\chi$ is the component of~$(e\setminus X)\cup\{\chi\}$ that contains~$\chi$ (in other words, the part of~$e$ that can be reached from~$\chi$ without passing through any other crossing). In this way an edge with~$\lambda$ crossings in~$\Gamma$ is assigned~$\lambda$ half-crossings. (Two halfcrossings of an edge~$e$ overlap iff the corresponding crossings are consecutive along~$e$.) The key observation is that at least a~$2/3$ fraction of all halfcrossings are incident to an endpoint of their edge. Only if an edge has three crossings, its halfsegment that corresponds to its middle crossing is not incident to an endpoint. 

  As we have~$2c$ halfcrossings in total, at most~$2c/3$ halfcrossings are not incident to an endpoint of their edge. Thus, for at least~$c-2c/3=c/3$ crossings both of its halfcrossings are incident to an endpoint. For any such crossing we can argue as in the proof of \cref{thm:cr-2planar} that by the maximality of~$\Gamma$ there is an uncrossed edge between the two endpoints, which are distinct because adjacent edges do not cross in~$\Gamma$. Every such edge can be obtained no more than twice, so that we find at least~$c/6$ uncrossed edges in~$\Gamma$, as claimed.
\end{proof}

\section{1-planar graphs}\label{sec:1planar}

In this section we focus on 1-planar graphs and we present lower and upper bounds on their edge density assuming that they are either $C_3$-free (\cref{ssec:1planar-c3free}) or $C_4$-free (\cref{ssec:1planar-c4free}) or of girth $5$ (\cref{ssec:1planar-girth5}).

\subsection{C${}_\mathbf{3}$-free 1-planar graphs}\label{ssec:1planar-c3free}

We start with the case of $C_3$-free 1-planar graphs, where we can derive an upper bound of $3(n-2)$ on their edge density (see \cref{thm:1-planar-C3-upper}); for a matching lower bound (up to a small additive constant) refer to \cite{DBLP:journals/dmgt/CzapPS16}.

\begin{theorem}\label{thm:1-planar-C3-upper}
Every $C_3$-free $1$-planar graph with~$n\ge 4$ vertices has at most~$3(n-2)$ edges.
\end{theorem}
\begin{proof}
We derive the upper bound by an application of the recently introduced edge-density formula for non-planar graphs~\cite{DBLP:journals/corr/abs-2311-06193}  given as follows:

\begin{equation}\label{eq:density-formula}    
|E| \leq t\left(|V|-2\right) - \sum_{c \in \mathcal{C}} \left(\frac{t-1}{4}||c||-t\right) - |\mathcal{X}|,
\end{equation}

\noindent where $\mathcal{C}$ and $\mathcal{X}$ denote the sets of cells and crossings, respectively. By setting $t=3$ to \eqref{eq:density-formula}, one gets $|E| \leq 3(n-2) + \frac{1}{2}|\mathcal{C}_5| - \frac{1}{2}|\mathcal{C}_6| - \ldots -  | \mathcal{X}|$, where $C_i$ denotes the set of cells of size $i$ with the size of a cell being the number of vertices and edge-segments on its boundary. Since each crossing is incident to at most two cells of size $5$ (as otherwise a $C_3$ is inevitably formed), it follows that $\frac{1}{2}|\mathcal{C}_5| \leq |\mathcal{X}|$, which  by the formula given above implies that $|E| \leq 3(n-2)$. 
\end{proof}

\subsection{C${}_\mathbf{4}$-free 1-planar graphs}\label{ssec:1planar-c4free}

We continue with the case of $C_4$-free 1-planar graphs. As in the case of $C_3$-free 1-planar graphs, we can again derive an upper bound of $3(n-2)$ for the edge-density using the density formula of \eqref{eq:density-formula}, since each crossing is incident to at most two cells of size $5$ (as otherwise a $C_4$ is formed). In the following theorem, we present an improved upper bound.

\begin{theorem}\label{thm:1-planar-C4-upper}
  Every $C_4$-free $1$-planar graph with~$n\ge 4$ vertices has at most~$\frac{5}{2}(n-2)$ edges.
\end{theorem}
\begin{proof}
We apply the discharging method with~$\alpha=4/5$ so that the statement follows by \eqref{eq:edge-bound}.
By \eqref{eq:remaining-charge} we have
\begin{equation}\label{eq:15bound}
\ch^-(f)=\frac15|\mathcal{V}(f)|+|f|-4\,.
\end{equation}

In particular, we have~$\ch^-(f)>0$ for all faces with at least four edge segments on the boundary. It remains to handle triangles. 

As the graph~$G$ is $1$-planar, every edge of~$G'$ is incident to at least one vertex in~$V$. It follows that
\begin{equation}\label{eq:vbound}
  |\mathcal{V}(f)|\ge\lceil|f|/2\rceil\,,
\end{equation}
for each~$f\in F'$. So every triangle~$f\in F'$ has either three vertices in~$V$ and~$\ch^-(f)=-2/5$ (type-1) or two vertices in~$V$ and one vertex 
in~$V'\setminus V$ with~$\ch^-(f)=-3/5$ (type-2).

We will argue how to make up for the deficits at triangles by transferring charges from neighboring faces.

First, let us discuss faces of size at least five. 
So consider~$f\in F'$ with~$|f|\ge 5$, and let~$k$ denote the number of triangles adjacent to~$f$ in the dual of~$G'$. 
Then for any vertex ~$v\in\mathcal{V'}(f)\setminus\mathcal{V}(f)$, at most one of the two edges incident to~$v$ along~$\partial f$ can be incident to a triangle of~$F'$ (because otherwise the two edges of~$G$ that cross at~$v$ induce a~$C_4$).
Thus,
\[
k\le|\mathcal{V}(f)|+\frac{|f|-|\mathcal{V}(f)|}{2}=
\frac{|f|+|\mathcal{V}(f)|}{2}\,.
\]

Together with \eqref{eq:15bound} we obtain
\[
\ch^-(f)=\frac15|\mathcal{V}(f)|+|f|-4=
\frac{|f|+|\mathcal{V}(f)|}{5}+\frac{4}{5}|f|-4\ge
\frac25 k\,,
\]
which shows that~$f$ can send a charge of~$2/5$ to every adjacent triangle.

\begin{figure}[htbp]
    \centering
    \begin{minipage}[b]{.3\textwidth}
    \centering  
    \includegraphics[page=3]{1-planar-cases.pdf}
    \subcaption{type-1}\label{fig:1-planar-cases:1}
    \end{minipage}
    \hfill
    \begin{minipage}[b]{.3\textwidth}
    \centering  
    \includegraphics[page=2]{1-planar-cases.pdf}
    \subcaption{type-2}\label{fig:1-planar-cases:2}
    \end{minipage}
    \hfill
    \begin{minipage}[b]{.3\textwidth}
    \centering  
    \includegraphics[page=1]{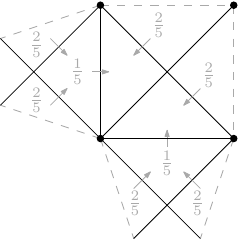}
    \subcaption{type-2}\label{fig:1-planar-cases:3}
    \end{minipage}
    \caption{Triangles in the planarization of~$C_4$-free $1$-planar graphs.}
    \label{fig:1-planar-cases}
\end{figure}

Next, consider a face~$f$ with~$|f|=4$. Combining \eqref{eq:remaining-charge} and \eqref{eq:vbound} we obtain~$\ch^-(f)=|\mathcal{V}(f)|/5\ge 2/5$. We claim that~$f$ can send a charge of~$2/5$ to every triangle that is adjacent to~$f$ via an edge of~$E'\setminus E$ and a charge of~$1/5$ to every triangle that is adjacent to~$f$ via an edge of~$E$. To see this, let us consider the three different types of quadrangles in~$F'$. By \eqref{eq:vbound} we have~$|\mathcal{V}(f)|\ge 2$. 

If~$|\mathcal{V}(f)|=2$, then there is at most one triangle adjacent to~$f$ because any two triangles adjacent to~$f$ induce a~$C_4$. So in this case~$f$ can send a charge of~$2/5$ to every adjacent~triangle. 

If~$|\mathcal{V}(f)|=3$, then any triangle adjacent to~$f$ via an edge of~$E'\setminus E$ induces a~$C_4$ in~$G$. Thus, there exist at most two triangles adjacent to~$f$ and every such triangle is adjacent via an edge of~$E$. So in this case~$f$ can send a charge of~$1/5$ to every adjacent triangle. 

Finally, if~$|\mathcal{V}(f)|=4$, then every triangle adjacent to~$f$ is adjacent via an edge of~$E$. As~$\ch^-(f)=\frac45$, also in this case~$f$ can send a charge of~$1/5$ to every adjacent triangle. This completes the proof of our claim.

So let us consider the incoming charges at triangles. For a type-1 triangle~$f$, neither of the adjacent faces is a type-1 triangle because such a pair would induce a~$C_4$ in~$G$. If at least two adjacent faces are type-2 triangles, then for each such triangle~$g$, neither of the other ($\ne f$) two faces adjacent to~$g$ are triangles because together with~$f$ and~$g$ they would induce a~$C_4$. It follows that~$g$ receives a charge of~$2\cdot 2/5=4/5$ from its two other ($\ne f$) neighbors, see \cref{fig:1-planar-cases:1}. As~$\ch^-(g)=-3/5$, the remaining charge of~$1/5$ can be passed on to~$f$. Then~$f$ receives a charge of~$2\cdot 1/5=2/5=-\ch^-(f)$ overall. 
Otherwise, at least two~of the three faces adjacent to~$f$ have size at least four. Each passes a charge of~$1/5$ across~the joint edge, which is in~$E$, to~$f$. So the deficit of~$\ch^-(f)=-2/5$ is covered in this~case as~well.

It remains to consider type-2 triangles. Let~$f$ be a type-2 triangle, and consider the two faces~$g_1,g_2$ that are adjacent to~$f$ via an edge of~$E'\setminus E$. If both~$g_1$ and~$g_2$ are triangles, then they induce a~$C_4$ in~$G$, in contradiction to~$G$ being~$C_4$-free. If both~$g_1$ and~$g_2$ have size at least four, then~$f$ receives a charge of~$2\cdot 2/5=4/5$ from them, which covers~$\ch^-(f)=-3/5$ and even leaves room to sent a charge of~$1/5$ across its third edge, which is in~$E$, see \cref{fig:1-planar-cases:2}. 

Hence, we may assume that without loss of generality~$g_1$ is a type-2 triangle and~$|g_2|\ge 4$. The third face~$g_3\notin\{g_1,g_2\}$ adjacent to~$f$ is not a type-1 triangle because then~$g_3$ together with~$g_1$ would induce a~$C_4$ in~$G$. If~$g_3$ is a type-2 triangle, then neither of its two other ($\ne f$) neighbors is a triangle because together with~$f$ and~$g_1$ there would be a~$C_4$ in~$G$. Therefore, we are in the case discussed above, where~$g_3$ receives a charge of~$4/5$ from its neighbors and passes on~$1/5$ to~$f$. Otherwise, we have~$|g_3|\ge 4$ and thus $g_3$ sends a charge of~$1/5$ to~$f$ across the joint edge, which is in~$E$. Together with the charge of~$2/5$ that~$f$ receives from~$g_2$ via the joint edge, which is in~$E'\setminus E$, this suffices to cover~$\ch^-(f)=-3/5$, see \cref{fig:1-planar-cases:3}. 
\end{proof}

\begin{restatable}{theorem}{oneplanarcfourlower}\label{thm:1planar-C4-lower}
    For every sufficiently large $n$, there exists a $C_4$-free $1$-planar graph on $n$ vertices with $2.4n-O(1)$ edges.
\end{restatable}

\begin{proof}
 Consider the following grid-based construction, see \cref{app:fig:1planar-C4-lower}. We put vertices at each coordinate $(i,j)$ for $i$ and $j$ integers. We call these vertices \emph{black} and depict them by black disks in \cref{app:fig:1planar-C4-lower}.
For $i$ odd and $j$ odd (that is, every second index), we put vertices at coordinates $(i+\frac{1}{2},j+\frac{1}{2})$.  We call these vertices \emph{red} and depict them by red squares in \cref{app:fig:1planar-C4-lower}.
For a graph with $n$ vertices this leads to $\frac{4}{5}n$ black vertices (at integer points) and $\frac{1}{5}n$ red vertices (at points at coordinates $(i+\frac{1}{2},j+\frac{1}{2})$). 

\begin{figure}[htbp]
    \hfill\nolinenumbers
    \begin{subfigure}[b]{.48\linewidth}
    \centering
    \includegraphics[page=4,scale=.8]{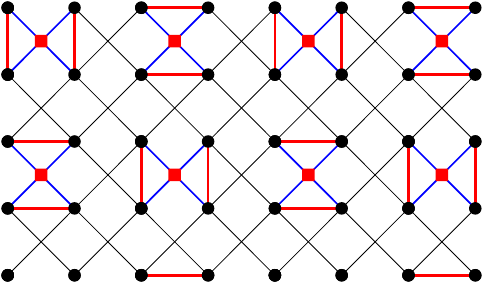}
    \subcaption{}
    \label{app:fig:1planar-C4-lower}
    \end{subfigure}
    \hfill
    \begin{subfigure}[b]{.48\linewidth}
    \centering
    \includegraphics[page=5,scale=.8]{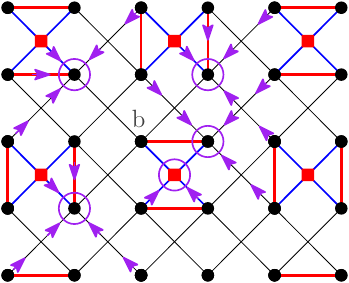}
    \subcaption{}
    \label{app:fig:1planar-C4-lower-c4free}
    \end{subfigure}
    \caption{(a)~A dense $C_4$-free $1$-plane graph.
    (b)~The neighborhood of a vertex $b$.}    
\end{figure} 

For each red vertex $r$ at position $(i+\frac{1}{2},j+\frac{1}{2})$, we add edges between $r$ and the vertex at position $(i,j)$, the vertex at position $(i+1,j)$, the vertex at position $(i+1,j+1)$, and the vertex at position $(i,j+1)$. (In Figure~\ref{app:fig:1planar-C4-lower}, these edges are depicted blue). 

For $i$ and $j$ both odd, we add an edge $(i,j)$ to $(i-1,j-1)$. For $i$ odd and $j$ even, we add an edge $(i,j)$ to $(i+1,j+1)$. For $i$ even, we add an edge $(i,j)$ to $(i-1,j-1)$ and an edge $(i,j)$ to $(i+1,j+1)$. (In Figure~\ref{app:fig:1planar-C4-lower}, these edges are depicted black).
(This construction so far can be seen as a grid consisting of crossings, where each second crossing is subdivided such that instead of a crossing point, there is a vertex of degree four.)  

 If $i$ is odd and $j \equiv 1 \;(\bmod\; 4)$ or $j \equiv 2 \;(\bmod\; 4)$, we add an edge from $(i,j)$ to $(i+1,j)$, and if $i$ is odd and $j \equiv 3 \;(\bmod\; 4)$, we add an edge from $(i,j)$ to $(i,j+1)$ and an edge from $(i+1,j)$ to $(i+1,j+1)$. (In Figure~\ref{app:fig:1planar-C4-lower}, these edges are depicted bold red).

We observe that this graph is symmetric in the sense that there are only two types of vertices. Choosing any black vertex, by rotating the drawing (by 90, 180 or 270 degrees) and mirroring it, it behaves like any other black vertex and the same holds for red vertices.

It is easy to see that this construction is 1-planar (and each crossing is the unique crossing between an edge connecting a vertex at position $(i,j)$ to one at position $(i+1,j+1)$ and an edge connecting a vertex at position $(i+1,j)$ to one at position $(i,j+1)$). To see that it is also $C_4$-free, we consider the $2$-neighborhoods of black vertices. Since every red vertex has only black vertices as neighbors, any $C_4$ contains at least two black vertices. Let $b$ be an arbitrary black vertex and let $N_b$ denote the set of the five neighbors of~$b$ (circled in \cref{app:fig:1planar-C4-lower-c4free}). Each~$v\in N_b$ has a set~$N_v$ of three or four neighbors other than~$b$. Any~$u\in N_v$ is either in~$N_b$ (forming a triangle~$bvu$) or~$v$ is the unique neighbor of~$u$ in~$N_b$ (indicated by arrows in \cref{app:fig:1planar-C4-lower-c4free}). Thus, there is no way to form a~$C_4$ passing through~$b$, and given that the choice of~$b$ was generic, the graph we built is~$C_4$-free.

Since each red vertex has degree four, and each black vertex has degree five, the number of edges is $(5\cdot\frac{4}{5}\cdot n+4\cdot\frac{1}{5}\cdot n)/2=2.4n$ minus some edges along the boundary. However, we can wrap the grid around a cylinder and make it constant width and arbitrary height, so as to obtain a constant size boundary.
\end{proof}

\subsection{1-planar graphs of girth 5}\label{ssec:1planar-girth5}

\begin{theorem}\label{thm:1-planar-girth-5}
  Every $1$-planar graph of girth~$5$ on~$n$ vertices has at most~$2.4n$ edges. 
\end{theorem}
\begin{proof}
We go through the proof of the 1-planar $C_4$-free case, and note that the arrow case as well as the type-1 triangles do not occur. If we choose $ \alpha = \frac{5}{6}$, the type-2 triangles have a negative charge of $-\frac{4}{6}$, and can get charges of $\frac{2}{6}$ from their immediate neighboring cells which are of  size at least 4. Note that for the case that those neighboring cells are of size 4, they have only one type-2 triangle by the $C_4$-freeness property, which suffices to provide enough charge. 
If a neighboring cell $c$ has size 5, then, by 1-planarity, it shares with at least one neighbor a planar edge, through which it does not have to contribute charge. Since the remaining charge of $c$ is at least $3 \cdot \frac{1}{6}+5-4 = 1.5$,
it can contribute to four neighbours $\frac{2}{6}$ charge each. 
If a neighboring cell $c$ is of size larger than 5, then its remaining charge is at least $\frac{1}{6}\vert \mathcal{V}(c)\vert  + \vert c \vert - 4 \ge \vert c \vert - 3 \ge \frac{2}{6} \vert c \vert$,
and therefore there is enough charge to provide $\frac{2}{6}$ charge to every neighboring type-2 triangle. 
This immediately gives that an $n$-vertex 1-planar graph of girth 5 has at most $\frac{12}{5}n = 2.4n$ edges.
\end{proof}

\begin{figure}[htb]
    \centering
    \includegraphics[page=4]{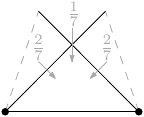}
    \caption{Triangles in the planarization of $1$-planar graphs of girth~$5$.}
    \label{fig:1-planar-girth-5}
\end{figure}

\begin{restatable}{theorem}{oneplanargirthfivelower}\label{thm:1planar-girth5-lower}    
For every sufficiently large $n$, there exists a $1$-planar graph of girth~$5$ on $n$ vertices with $(2 + \frac{1}{6})n-O(1)\approx 2.167n-O(1)$ edges.
\end{restatable}    

\begin{proof}
The construction is illustrated in \cref{app:fig:1planar-girth5-lower:1}. If we ignore the red edges, then we have a hexagonal grid each containing a uniform pattern of 4 vertices and 9 edges. In particular, each tile (incl.~the boundary) of this grid induces a $1$-plane drawing of the Petersen graph. 
Let $n'$ denote the number of vertices spanned by the hexagonal grid. Then, by Euler's formula, we have $\frac{n'}{2}$ hexagons. Together with the $\frac{n'}{2}$ disjoint patterns, we
have $3n'$ vertices and $\frac{3}{2}\cdot n' + 9\cdot  \frac{n'}{2} = 6n'$ edges. Assume that $n= 3n'$. Then, the obtained graph has $n$ vertices and $2n$ edges distributed among $\frac{n}{6}$ hexagons.
Additionally, for every triple of hexagons, there exists one red edge that joins two of these hexagons and four other red edges that lead to hexagons of other triples of hexagons. Amortized over all the hexagons, we can account one additional red edge per hexagon. This gives~$\frac{n}{6}$ additional red edges, for a total of~$\left(2+\frac{1}{6}\right)n$ edges. We argue as in the proof of \cref{thm:1planar-C4-lower} that the boundary effects result in a loss of a constant number of edges only.

 \begin{figure}[htbp]
     \begin{subfigure}[b]{0.65\linewidth}
    \center
    \includegraphics[width=\textwidth, page=2]{figures/1-planar-girth-5-2.166n.pdf}
    \subcaption{}
    \label{app:fig:1planar-girth5-lower:1}
    \end{subfigure} 
    \hfill
    \begin{subfigure}[b]{0.3\linewidth}
    \center
    \includegraphics[width=\textwidth, page=1]{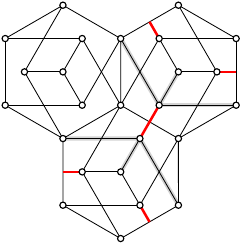}
    \subcaption{}
    \label{app:fig:1planar-girth5-lower:2}
    \end{subfigure}  
    \hfill
    \caption{(a) A $1$-planar graph of girth~$5$ with about $2.167n$ edges (\cref{thm:1planar-girth5-lower}). The construction consists of repeated triplets of hexagonal tiles (bordered by thick edges, also shown in (b)).}
    \label{app:fig:1planar-girth5-lower}
\end{figure}

It remains to argue that the constructed graph has girth~$5$. First, observe that the subgraph within each hexagonal tile has girth~$5$ because it is a Petersen graph, which is known to have girth~$5$. It follows that every~$C_3$ or~$C_4$, if any, uses vertices from at least two different tiles.
Second, we argue that no~$C_3$ or~$C_4$ uses a red edge. To see this consider the neighbors of the two endpoints of a red edge and observe that they are at pairwise distance at least two; see \cref{app:fig:1planar-girth5-lower:2}, where the edges in the neighborhood of a red edge are shown gray. As red edges are the only edges that cross tile boundaries and boundary edges are shared among adjacent tiles, it follows that every~$C_3$ or~$C_4$, if any, uses at least two nonadjacent vertices~$u,v$ on the boundary~$\partial T$ of a tile~$T$ and exactly one vertex~$z$ in the interior of $T$. Then~$u$ and~$v$ are antipodal on~$\partial T$ (i.e., at distance three along~$\partial T$). In particular, these vertices do not form a~$C_3$. Further, the tile~$T$ is the unique common tile of~$u$ and~$v$, so there is no common neighbor of~$u$ and~$v$ outside of~$T$. As~$z$ is the only common neighbor of~$u$ and~$v$ inside~$T$, it follows that there is no~$C_4$ through~$u,v,z$. Thus, our graph has girth~$5$.
\end{proof}

\section{2-planar graphs}\label{sec:2planar}

In this section, we focus on $2$-planar graphs and we present bounds on their edge density assuming that they are $C_3$-free (\cref{ssec:2planar-c3free}) or $C_4$-free (\cref{ssec:2planar-c4free}) or of girth~$5$ (\cref{ssec:2planar-girth5}). 

\subsection{C${}_\mathbf{3}$-free 2-planar graphs} \label{ssec:2planar-c3free}

For the maximum edge density of $C_3$-free 2-planar graphs, we can derive an upper bound of $4(n-2)$ (see \cref{thm:2-planar-C3-upper}); for a lower bound of $3.5(n-2)$ refer to \cite{DBLP:conf/isaac/AngeliniB0PU18}.

\begin{theorem}\label{thm:2-planar-C3-upper}
$C_3$-free $2$-planar graphs with $n$ vertices have at most $4(n-2)$ edges.
\end{theorem}
\begin{proof}
To derive the upper bound, we apply the discharging method with~$\alpha=\frac{1}{2}$ so that the statement follows by \eqref{eq:edge-bound}. By \eqref{eq:remaining-charge} we have
\begin{equation}
\ch^-(f)=\frac12|\mathcal{V}(f)|+|f|-4\,.
\end{equation}

In particular, we have~$\ch^-(f) \ge 0$ for all faces with at least four edges on the boundary. It remains to handle triangles. Since we consider $C_3$-free graphs, we distinguish between three types of triangles; those with 0, 1 and 2 vertices on their boundaries and it is not difficult to observe that the latter ones have zero charge, while the former ones have charge $-1$ and $-\frac12$, respectively.

\begin{figure}[htbp]
    \centering
    \begin{subfigure}[b]{0.3\textwidth}
    \centering  
    \includegraphics[page=1]{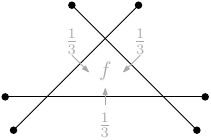}
    \subcaption{}
    \label{fig:2-planar-cases:1}
    \end{subfigure}
    \hfill
    \begin{subfigure}[b]{0.3\textwidth}
    \centering  
    \includegraphics[page=2]{2-planar-cases.pdf}
    \subcaption{}\label{fig:2-planar-cases:2}
    \end{subfigure}
    \hfill
    \begin{subfigure}[b]{0.3\textwidth}
    \centering  
    \includegraphics[page=3]{2-planar-cases.pdf}
    \subcaption{}\label{fig:2-planar-cases:3}
    \end{subfigure}
    
    \begin{subfigure}[b]{0.3\textwidth}
    \centering  
    \includegraphics[page=4]{2-planar-cases.pdf}
    \subcaption{}\label{fig:2-planar-cases:4}
    \end{subfigure}
    \begin{subfigure}[b]{0.3\textwidth}
    \centering  
    \includegraphics[page=5]{2-planar-cases.pdf}
    \subcaption{}\label{fig:2-planar-cases:5}
    \end{subfigure}
    \caption{Triangles in the planarization of~$C_3$-free $2$-planar graphs}
    \label{fig:2-planar-cases}
\end{figure}

For each triangle $f$ with zero or one vertices on its boundary, our strategy is to transfer at least $\frac{1}{4}$ and at most $\frac{1}{3}$ units of charge from the cells neighboring $f$. 
If $f$ has no vertices on its boundary, then we will transfer $\frac{1}{3}$ units of charge from each neighboring cell. Otherwise, we will transfer $\frac{1}{4}$ units of charge each from two neighboring cells of $f$; see \cref{fig:2-planar-cases}. 

Assume first that $\vert \mathcal V(f) \vert = 0$; see~\cref{fig:2-planar-cases:1}. Since $f$ is triangular, it follows that $f$ is formed by three mutually crossing edges. Our strategy is to transfer $\frac{1}{3}$ units of charge from each cell neighboring $f$. Since $f$ neighbors three such cells, this is enough to bring the remaining charge of $f$ from $-1$ to $0$. Let $c$ be a neighboring cell of $f$. It follows that $\vert c \vert \geq 4$ with two vertices on its boundary. 
Thus its remaining charge is: 
\[ \frac{1}{2}\vert \mathcal{V}(c)\vert  + \vert c \vert - 4 \ge 1 -4 + \vert c \vert  = \vert c \vert - 3 
\]
This implies that if $\vert c \vert \geq 5$, then the remaining charge of $c$ is at least $2$, in which case $c$ can transfer $\frac{1}{3}$ units of charge to $f$ and its remaining charge will be enough to distributed to the rest of its neighboring cells.
For the second case, we assume $\vert c \vert = 4$. Since $c$ has two vertices that appear consecutively on its boundary, it follows that one of the sides that bound $c$ is crossing free. Denote this side by $e$ and let $c'$ be the cell on the other side of $e$. It follows that $c'$ is neither a triangle with zero vertices nor a triangle with one vertex on its boundary. Hence, there is no need to transfer charge from $c$ to $c'$ according to our strategy. It follows that there are at most $3$ neighboring cells that $c$ may have to transfer charge to. Hence, $c$ can transfer $\frac{1}{3}$ units of charge to $f$ and its remaining charge will be enough to distributed to the rest of its neighboring cells, if needed.

To complete the proof of the theorem, we next consider the case in which $|\mathcal{V}(f)|=1$; see~\cref{fig:2-planar-cases:2}. Let $u$ be the vertex on the boundary of $f$ and let $(u_1,u_2)$ be the edge with one of its segments on the boundary of $f$. Let $c_1$ and $c_2$ be the two neighboring cells of $f$ that share the two sides of $f$ incident to its vertex. Since we consider $C_3$-free graphs, it follows that $(u,u_1)$ and $(u,u_2)$ cannot be both in the graph. Assume that $(u,u_i)$ with $i \in \{1,2\}$ is not part of the graph. Then, the corresponding cell $c \in \{c_1,c_2\}$ neighboring $f$ and having vertex $u$ and $u_i$ on its boundary has size at least $4$, which means that its remaining charge is at least:  
\[ \frac{1}{2}\vert \mathcal{V}(c)\vert  + \vert c \vert - 4 \ge 1 -4 + \vert c \vert  = \vert c \vert - 3 \ge \frac{1}{4} \vert c \vert
\]
Hence, we can safely transfer $\frac{1}{4}$ units of charge from $c$ to $f$, since the remaining charge of $c$ would be enough for being distributed to the remaining cells neighboring $c$, if needed. This implies that if both $(u,u_1)$ and $(u,u_2)$ are not in the graph, then each of the cells $c_1$ and $c_2$ can transfer $\frac{1}{4}$ units of charge to $f$ and then we are done. So, in the rest we can assume that this is not the case. 

Let $c'$ be the face neighboring $f$ that is on the other side of the edge $(u_1,u_2)$. If $c'$ has at least two vertices on its boundary, then as above we transfer $\frac{1}{4}$ units of charge from $c'$ to $f$ and the remaining charge of $c'$ would be enough for being distributed to the remaining cells neighboring $c'$, if needed. So, it remains to consider the cases in which $c'$ has either no or one vertex on its boundary. 

Assume first that $c'$ has one vertex on its boundary, that is, $\mathcal V(c')=1$. Then: 
\[ \frac{1}{2}\vert \mathcal{V}(c')\vert  + \vert c' \vert - 4 \ge \frac{1}{2} -4 + \vert c' \vert  = \vert c' \vert - 3.5 
\]
If $c'$ is such that $\vert c' \vert \geq 5$, then $\vert c' \vert - 3.5 \ge \frac{1}{4} \vert c' \vert$ holds and as above we can safely transfer $\frac{1}{4}$ units of charge from $c'$ to $f$. So, it remains to argue for the cases in which $\vert c' \vert \in \{3,4\}$. First, we observe that $\vert c' \vert \neq 3$, as otherwise the two edges incident to $u$ bounding $f$ would form a pair of parallel edges. Hence, we may assume that $\vert c' \vert = 4$; see \cref{fig:2-planar-cases:3}. Since $c'$ has one vertex on its boundary,  its remaining charge is $\frac{1}{2}$. In this case, we argue that at most two neighboring cells, namely, $f$ and another one, may need additional charge from $c'$. In particular, the two cells neighboring $c'$ that have the vertex of $c'$ on their boundary do not need additional charge, since none of them can be a triangle with zero or one vertex on its boundary. This means that we can safely transfer $\frac{1}{4}$ units of charge from $c'$ to $f$, as desired. 

To complete the case analysis, we need to consider the case that $c'$ has no vertex on its boundary. In this case, the remaining charge of $c'$ is $\vert c' \vert - 4$. 
If $\vert c' \vert \geq 6$, then the remaining charge of $c'$ is at least $2$, which implies that $\frac{1}{4}$ units of charge can be safely transferred to $f$ and the remaining charge of $c'$ will be enough for being distributed to the rest of the cells neighboring $c'$, if needed. So, we may assume that $\vert c' \vert \in \{3,4,5\}$. First, we observe that $\vert c' \vert \neq 3$, as otherwise the two edges incident to $u$ bounding $f$ would form a pair of crossing edges, which is not possible in simple drawings. Hence, $\vert c' \vert \in \{4,5\}$. If $\vert c' \vert = 4$, then its remaining charge is $0$ and clearly it cannot transfer charge to $f$. In this case, we consider the cell $c''$ neighboring $c'$, which does not share a crossing point with $f$; see~\cref{fig:2-planar-cases:4}. It follows that $\vert c'' \vert \ge 4$ and $c''$ has two vertices on its boundary. Since $c'$ does not require a transfer of charge, we transfer $\frac{1}{4}$ units of charge from $c''$ to $f$ and as in the first case of the proof the remaining charge of $c''$ for being distributed to the rest of the cells neighboring $c''$.  

To complete the proof of the case $\vert \mathcal V(f) \vert = 1$, consider now the case $\vert c' \vert = 5$. In this case, the remaining charge of $c'$ is $1$ and this is enough to contribute a $\frac{1}{4}$ to at most four neighboring cells. Hence, we may assume that $c'$ has to transfer $\frac{1}{4}$ units of charge to exactly five neighboring cells; see~\cref{fig:2-planar-cases:5}. In this case, it follows that none of the edges $(u,u_1)$ and $(u, u_2)$ is part of the graph (as otherwise there is a $C_3$; a contradiction). However, we have assumed that one of these edges belongs to the graph. 
\end{proof}

\subsection{C${}_\mathbf{4}$-free 2-planar graphs} \label{ssec:2planar-c4free}

We continue with the case of $C_4$-free 2-planar graphs, deriving an upper bound $3.929n$ on their maximum edge density (\cref{thm:2planar-c4}); for a lower bound of $2.5n -O(1)$ refer to~\cref{thm:2planar-c4-lower}.

\begin{theorem}\label{thm:2planar-c4} 
  Every $C_4$-free $2$-planar graph on~$n\ge 2$ vertices has at most
  \[
  \sqrt[3]{\frac{190,125}{3,136}}n<3.929n
  \]
  edges.
\end{theorem}
\begin{proof}
  Let~$G$ be a $C_4$-free graph with~$n$ vertices and~$m$ edges, and let~$\Gamma$ be a minimum-crossing drawing of~$G$. Then \cref{thm:naivelowerbound} in combination with the upper bound of $\frac{15}{7}(n-2)$ by Dowden~\cite{DBLP:journals/jgt/Dowden16} regarding the edge density of $C_4$-free planar graphs and  
  \cref{thm:1-planar-C4-upper} yields:
  \[
  \mathrm{cr}(G)\ge 2m-\frac{5}{2}n-\frac{15}{7}n=
2m-\frac{65}{14}n\,.
  \]
  By applying \cref{thm:crossing-lemma} for $a=2$ and $b=\frac{65}{14}$, we obtain the following lower bound on the number of crossings of $G$ when  $m\ge \frac{195n}{56}\approx 3.482n$.
  \begin{equation}\label{eq:crossing-lemma-c4-free}
  \mathrm{cr}(G)\ge \frac{6,272 m^3}{114,075n^2}\,.
  \end{equation}
  Assume now that $G$ is additionally $2$-planar. Then by \cref{thm:cr-2planar}, we obtain $\mathrm{cr}(G) \leq \frac{10n}{3}$. Hence, by \eqref{eq:crossing-lemma-c4-free} we have:
  \[
  \frac{6,272 m^3}{114,075n^2}\le\frac{10}{3}n\iff
  m^3\le\frac{190,125}{3,136}n^3\,.\qedhere
  \]
\end{proof}

\begin{corollary}\label{cor:c4-free-k-planar-1}
Every $C_4$-free $k$-planar graph on~$n\ge 2$ vertices and $m \geq 3.483n$ edges has at most
  \[
  \sqrt{\frac{114,075}{12,544}} \cdot \sqrt{k} \cdot n<
  3.016\sqrt{k}n
  \]
  edges.
\end{corollary}
\begin{proof}
    Let $G$ be a $C_4$-free $k$-planar graph with $n$ vertices and  $m \geq 3.483n$ edges. 
    By \eqref{eq:crossing-lemma-c4-free}, we know a lower bound on its number of crossings, namely, 
    \[
   \mathrm{cr}(G)\ge \frac{6,272 m^3}{114,075n^2}\,.
   \]
   On the other hand, since $G$ is $k$-planar, it holds $\frac{km}{2} \geq cr(G)$. Combining those, we get:
   \[
  \frac{12,544 m^2}{114,075n^2}\le k\iff
  m\le\sqrt{\frac{114,075}{12,544}k}n\,.\qedhere
  \]
\end{proof}

\begin{remark}\label{thm:c4-free-k-planar-2}    
An asymptotically better bound of~$\Theta(\sqrt[3]{k}n)$ edges, 
which however holds for significantly denser graphs only, can be obtained by combining an improved crossing lemma for $C_4$-free graphs  by Pach, Spencer, and T{\'{o}}th~\cite[Theorem 3.1]{pach1999new} with the trivial upper bound of at most~$km/2$ crossings for~$k$-planar graphs.
\end{remark}

\begin{restatable}{theorem}{twoplanarcfourlower}\label{thm:2planar-c4-lower}     For every sufficiently large~$n$, there exists a $C_4$-free $2$-planar graph on~$n$ vertices with~$2.5n-O(1)$ edges.
\end{restatable}

\begin{figure}[htbp]
  \centering
  \includegraphics[page=1]{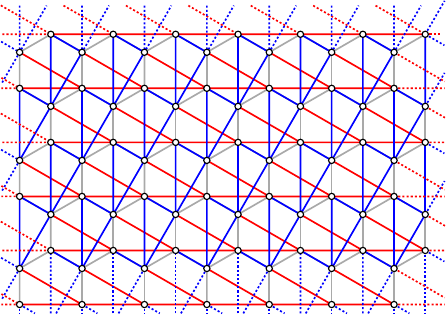}
  \includegraphics[page=2]{2-planar-c4-free-2.5n.pdf}
  \caption{A $2$-plane graph 
  with~$\approx 2.5n$ edges, shown red and blue. Gray shows the grid only.}
  \label{app:fig:2planar-c4-lower}
\end{figure}

\begin{proof}
  We arrange the vertices so that they form a hexagonal grid and connect the vertices in each grid cell by five edges; see \cref{app:fig:2planar-c4-lower}. More precisely, if we denote the vertices on the boundary of a grid cell in anticlockwise order, starting with a fixed direction by~$1,2,3,4,5,6$, then we add the edges~$14,24,35,36,56$. 
  It is easily checked that the resulting drawing is~$2$-plane, and that every internal (sufficiently far away from the boundary) vertex of the grid has degree five. Thus, the number of edges is about~$5n/2$. We argue as in the proof of \cref{thm:1planar-C4-lower} that the boundary effects result in a loss of a constant number of edges only.
  
  It remains to argue that the graph is~$C_4$-free. To see this, first note that the graph is highly symmetric, and we have only two types of vertices: \emph{blue} vertices are incident to two triangles and \emph{red} vertices are incident to one triangle only. Both classes of vertices in isolation induce a collection of disjoint paths. So any~$C_4$ has to use at least one blue vertex; let us denote this vertex by~$v$. Consider the $2$-neighborhood of~$v$: Let~$N_v$ denote the set of the five neighbors of~$v$ (circled in \cref{app:fig:2planar-c4-lower}). Each~$u\in N_v$ has a set~$N_u$ of four neighbors other than~$v$. Any~$w\in N_u$ is either in~$N_v$ (forming a triangle~$uvw$) or~$u$ is the unique neighbor of~$w$ in~$N_v$ (indicated by an arrow in \cref{app:fig:2planar-c4-lower}). Thus, there is no way to form a~$C_4$ passing through~$v$, and given that the choice of~$v$ was generic, the graph we built is~$C_4$-free.
\end{proof}

\subsection{2-planar graphs of girth 5}\label{ssec:2planar-girth5}

We conclude \cref{sec:2planar} with the case of 2-planar graphs of girth $5$.

\begin{theorem}\label{thm:2planar-g5} 
  Every $2$-planar graph of girth~$5$ on~$n$ vertices has at most
  \[
  \sqrt[3]{\frac{11,163}{240}}n<3.597n
  \]
  edges.
\end{theorem}

\begin{proof}
  Let~$G$ be a graph of girth~$5$ with~$n$ vertices and~$m$ edges. As a consequence of Euler's Formula, every planar graph of girth~$g$ on~$n\ge 3$ vertices has at most~$g(n-2)/(g-2)$ edges. Plugging this together with \cref{thm:1-planar-girth-5} into \cref{thm:naivelowerbound} we get
  \[
  \mathrm{cr}(G)\ge 2m-\frac{5}{3}n-\frac{12}{5}n=2m-\frac{61}{15}n\,.
  \]
  By applying \cref{thm:crossing-lemma} for $a=2$ and $b=\frac{61}{15}$, we obtain the following lower bound on the number of crossings of $G$ when  $m\ge \frac{61n}{20}$.
  \begin{equation}\label{eq:2-planar-girth-5}
  \mathrm{cr}(G)\ge\frac{800 m^3}{11,163 n^2}\,.
  \end{equation}
  Assume now that $G$ is additionally $2$-planar. Then by \cref{thm:cr-2planar}, we obtain $\mathrm{cr}(G) \leq \frac{10n}{3}$. Hence, by \eqref{eq:2-planar-girth-5} we have:
  \[
  \frac{800 m^3}{11,163 n^2}\le\frac{10}{3}n\iff
  m^3\le\frac{11,163}{240}n^3\,.\qedhere
  \]
\end{proof}

\noindent The next corollary follows from \eqref{eq:2-planar-girth-5} of the proof of \cref{thm:2planar-g5}; its proof is analogous to the one of \cref{cor:c4-free-k-planar-1}.

\begin{restatable}{corollary}{gfivekplanarone}\label{cor:g5-k-planar-1}
Every $k$-planar graph of girth $5$ on~$n\ge 2$ vertices and $m \geq 3.05n$ edges has at most
  \[
  \sqrt{\frac{11,163}{1,600}} \cdot \sqrt{k} \cdot n<
  2.642\sqrt{k}n
  \]
  edges.
\end{restatable}
\begin{proof}
    Let $G$ be a $k$-planar graph of girth $5$ with $n$ vertices and  $m \geq 3.05n$ edges. 
    By \eqref{eq:2-planar-girth-5}, we know a lower bound on its number of crossings, namely, 
    \[
   \mathrm{cr}(G)\ge \frac{800 m^3}{11,163n^2}\,.
   \]
   On the other hand, since $G$ is $k$-planar, it holds $\frac{km}{2} \geq \mathrm{cr}(G)$. Combining those, we get:
   \[
  \frac{1,600 m^2}{11,163n^2}\le k\iff
  m\le\sqrt{\frac{11,163}{1,600}k}n\,.\qedhere
  \]
\end{proof}

\begin{restatable}{theorem}{twoplanargfivelower}\label{thm:2planar-g5-lower}
    For every sufficiently large $n$, there exists a $2$-planar graph of girth~$5$ on $n$ vertices with $(2 + \frac{2}{7})n-O(1)\approx 2.286n-O(1)$ edges.
\end{restatable}
\begin{proof}
The construction is illustrated in \cref{app:fig:2planar-girth5-lower:1}. Starting with the construction from the proof of \cref{thm:1planar-girth5-lower}, we add three blue vertices in every third hexagon and three edges each starting from a blue vertex and ending in different hexagonal tiles. 
Let $n'$ denote the number of vertices spanned by the hexagonal grid. Then, by Euler's formula, we have $\frac{n'}{2}$ hexagons. Ignoring the red and blue edges, we have $3n'$ vertices and $6n'$ edges. As pointed out in the proof of \cref{thm:1planar-girth5-lower}, there is amortized one additional red edge per hexagon. With an analogue analysis, we obtain that there are three blue vertices and nine blue edges per a triple of hexagons.
This gives in total $3n' + \frac{n'}{2}$ vertices and $6n' + \frac{4n'}{2}$ edges. Setting $n=3.5n'$, we obtain that the graph has $\frac{16n}{7}$ edges.
We argue as in the proof of \cref{thm:1planar-C4-lower} that the boundary effects the result in a loss of a constant number of edges only.

 \begin{figure}[htbp]
    \begin{subfigure}[b]{0.65\linewidth}
    \center
    \includegraphics[width=\textwidth, page=2]{figures/2-planar-girth-5-2.2857n.pdf}
    \subcaption{}
    \label{app:fig:2planar-girth5-lower:1}
    \end{subfigure} 
    \hfill
    \begin{subfigure}[b]{0.3\linewidth}
    \center
    \includegraphics[width=\textwidth, page=1]{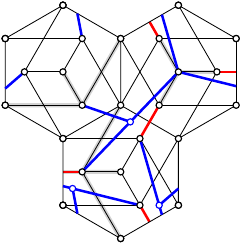}
    \subcaption{}
    \label{app:fig:2planar-girth5-lower:2}
    \end{subfigure}  
    \hfill
    \caption{Illustration for the proof of \cref{thm:2planar-g5-lower}. The construction consists of repeated triplets of hexagonal tiles (bordered by thick edges, also shown in (b)).}
    \label{app:fig:2planar-girth5-lower}
\end{figure}
It remains to show that the graph has girth 5. Assume $C$ is a cycle of length three or four. By the proof of \cref{thm:1planar-girth5-lower}, $C$ includes a blue edge $uv$ with $v$ is a blue vertex. As only blue edges end at a blue vertex, there exist a second blue edge $vw$ in $C$. We consider the neighbors of $u$ and $w$ and observe that they are at pairwise distance at least one; see \cref{app:fig:2planar-girth5-lower:2}, where the edges in the neighborhood of blue edges are shown gray.
\end{proof}

\section{3-planar graphs}\label{sec:3planar}

This section is devoted to $3$-planar graphs and is structured analogously to \cref{sec:2planar}. 

\subsection{C${}_\mathbf{3}$-free 3-planar graphs} \label{ssec:3planar-c3free}

\begin{theorem}\label{thm:3planar-c3}
  Every $C_3$-free $3$-planar graph on~$n\ge 2$ vertices has at most
  \[  
  \sqrt[3]{\frac{2,673}{20}}n<5.113n
  \]
  edges.
\end{theorem}
\begin{proof}
  Let~$G$ be a $C_3$-free 
  graph with~$n$ vertices and~$m$ edges, and let~$\Gamma$ be a minimum-crossing drawing of~$G$. Then \cref{thm:naivelowerbound} in combination with the fact that $C_3$-free $n$-vertex planar graphs with at most $2n-4$ edges and \cref{thm:1-planar-C4-upper} yields
  \[
  \mathrm{cr}(G)\ge 3m-2n-3n-4n=
3m-9n\,.
  \]
  By applying \cref{thm:crossing-lemma} for $a=3$ and $b=9$, we obtain the following lower bound on the number of crossings of $G$ when  $m\ge \frac{9n}{2}$.
  \begin{equation}\label{eq:crossing-lemma-c3-free}
  \mathrm{cr}(G)\ge \frac{4 m^3}{81 n^2}\,.
  \end{equation}
  Assume now that $G$ is additionally $3$-planar. Then by \cref{thm:cr-3planar}, we obtain $\mathrm{cr}(G) \leq \frac{33n}{5}$. Hence, by \eqref{eq:crossing-lemma-c3-free} we have:
  \[
  \frac{4 m^3}{81n^2}\le\frac{33}{5}n\iff
  m^3\le\frac{2,673}{20}n^3\,.\qedhere
  \]
\end{proof}

\begin{corollary}\label{cor:c3-k-planar-1}
Every $C_3$-free $k$-planar graph on~$n\ge 2$ vertices and $m \geq \frac{9}{2}n$ edges has at most
  \[
  \sqrt{\frac{81}{8}} \cdot \sqrt{k} \cdot n<
  3.182\sqrt{k}n
  \]
  edges.
\end{corollary}
\begin{proof}
    Let $G$ be a $C_3$-free $k$-planar graph with $n$ vertices and  $m \geq \frac{9}{2}n$ edges. 
    By \eqref{eq:crossing-lemma-c3-free}, we know a lower bound on its number of crossings, namely, 
    \[
   \mathrm{cr}(G)\ge \frac{4 m^3}{81n^2}\,.
   \]
   On the other hand, since $G$ is $k$-planar, it holds $\frac{km}{2} \geq \mathrm{cr}(G)$. Combining those, we get:
   \[
  \frac{8 m^2}{81n^2}\le k\iff
  m\le\sqrt{\frac{81}{8}k}n\,.\qedhere
  \]
\end{proof}

\subsection{C${}_\mathbf{4}$-free 3-planar graphs} \label{ssec:3planar-c4free}

\begin{theorem}\label{thm:3planar-c4-free}
  Every $C_4$-free $3$-planar graph on~$n\ge 2$ vertices has at most
  \[
  \sqrt[3]{\frac{3,764,475}{31,360}}n<4.933n
  \]
  edges.
\end{theorem}
\begin{proof}
  Let~$G$ be a $C_4$-free 
  graph with~$n$ vertices and~$m$ edges. By \eqref{eq:crossing-lemma-c4-free}, we have: 
  \begin{equation}\label{eq:crossing-lemma-c4-free-repeat}
  \mathrm{cr}(G)\ge \frac{6,272 m^3}{114,075n^2}\,.
  \end{equation}
  Assume now that $G$ is additionally $3$-planar. Then by \cref{thm:cr-3planar}, we obtain $\mathrm{cr}(G) \leq \frac{33n}{5}$. Hence, by \eqref{eq:crossing-lemma-c4-free-repeat} we have:
  \[
  \frac{6,272 m^3}{114,075n^2}\le\frac{33}{5}n\iff
  m^3\le\frac{3,764,475}{31,360}n^3\,.\qedhere
  \]
\end{proof}

\begin{theorem}\label{thm:3planar-g5}
  Every $3$-planar graph of girth~$5$ on~$n$ vertices has at most
  \[
  \sqrt[3]{\frac{368,379}{4,000}}n < 4.516n
  \]
  edges.
\end{theorem}
\begin{proof}
  Let~$G$ be a  
  graph of girth 5 with~$n$ vertices and~$m$ edges. By \eqref{eq:2-planar-girth-5}, we have: 
  \begin{equation}\label{eq:3-planar-girth-5-repeat}
  \mathrm{cr}(G)\ge \frac{800 m^3}{11,163n^2}\,.
  \end{equation}
  Assume now that $G$ is additionally $3$-planar. Then by \cref{thm:cr-3planar}, we apply the upper bound for the crossing number $\mathrm{cr}(G) \leq \frac{33n}{5}$. Hence, by \eqref{eq:3-planar-girth-5-repeat} we have:
  \[
  \frac{800 m^3}{11,163 n^2}\le\frac{33}{5}n\iff
  m^3\le\frac{368379}{4000}n^3\,.\qedhere
  \]
\end{proof}

\begin{theorem}\label{thm:3planar-g5-lower}
    For every sufficiently large $n$, there are $3$-planar graphs of girth $5$ with $2.5n - O(1)$ edges.
\end{theorem}
\begin{proof}
We give a construction achieving that density, see \cref{app:fig:3planar-girth5-lower_together}(a).
Consider the plane like a grid. We put vertices on all integer points $(i,j)$. We call vertices at positions $(i,j)$ where $i$ and $j$ have the same parity, \emph{red vertices}, and the other vertices \emph{black vertices}.  

\begin{figure}[htbp]
    \centering
    \includegraphics[page=3]{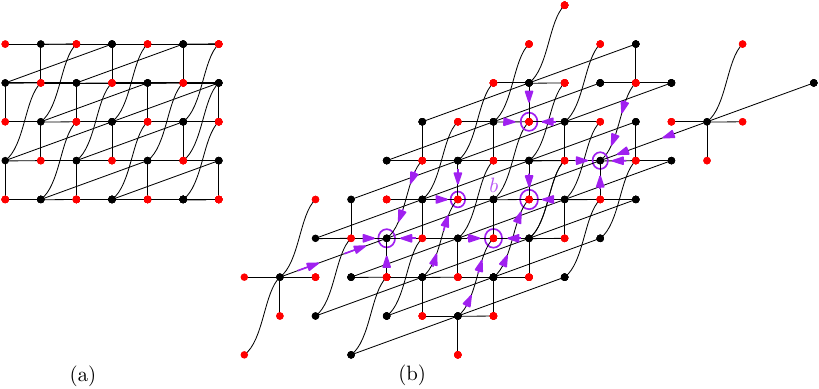}
    \caption{3-planar $C_4$-free construction. (a) Overview of the construction. (b) The vertex $b$, its neighbors (circled), vertices that are the neighbors of $N(b)$ and its neighbors. Each neighbor of $N(b)$ has a unique neighbor in $N(b)$ indicated by arrows.}
    \label{app:fig:3planar-girth5-lower_together}
\end{figure}

We put edges always between a vertex on point $(i,j)$ and a vertex on point $(i,j+1)$. 
We connect each black vertex at position $(i,j)$ via an edge with the (red) vertex at position $(i+2,j+1)$ and via an edge with the (black) vertex at position $(i+1, j+3)$. 
Finally, we connect each red vertex at position $(i,j)$ via an edge with the (black) vertex from $(i,j)$ to $(i+1,j)$. 

This way, all $\frac{n}{2}$ vertices at positions $(i,j)$ where the parity of $i$ and $j$ is different have degree six, while those at position $(i,j)$ where the parity of $i$ and $j$ is the same have degree four.  In total, this gives a density of $\frac{6}{4}+\frac{4}{4}=2.5$.

It can be observed that the construction is $3$-planar.

To see that the drawing is of girth five, we first observe that the drawing is symmetric in the sense that there are only two classes of vertices that all red vertices behave the same as each other, and all black vertices behave the same as each other. The four neighbors of a red vertex are all black. Thus, each cycle contains at least two black vertices. Let $b$ be an arbitrary black vertex and let $N_b$ denote the set of the six neighbors of~$b$ (circled in \cref{app:fig:3planar-girth5-lower_together}(b)). Each~$v\in N_b$ has a set~$N_v$ of three or five neighbors other than~$b$. No~$u\in N_v$ is in~$N_b$, and thus there is no way to form a~$C_3$ passing through~$b$. Further, $v$ is the unique neighbor of~$u$ in~$N_b$ (indicated by arrows in \cref{app:fig:3planar-girth5-lower_together}(b)). Thus, there is no way to form a~$C_4$ passing through~$b$. Given that the choice of~$b$ was generic, the graph we built is~$C_3$-free and $C_4$-free and thus has girth five.
\end{proof}

\section{Conclusions and open problems}
\label{sec:conclusions}

In this work, we continued an active research branch in Graph Drawing seeking for new edge density bounds for $k$-planar graphs that avoid certain forbidden substructures, namely, cycles of length $3$ or $4$ or both of them. For each of these settings, our focus was on $k$-planar graphs, with $k \in \{1,2,3\}$, as well as on general $k$. Several open problems have been triggered:

\begin{itemize}
\item The first one is the obvious one, that is, to close the gaps between the lower and the upper bounds reported in \cref{tab:bounds}. We believe that this is a challenging open problem. 
\item In particular, it seems to us that the lower bounds for $2$- and $3$-planar can be improved. 
\item Note that there is a lot of empty space to fill in \cref{tab:crossing-lemma} where we did not find any reasonably good bounds.
\item Another promising research direction is to study the edge density of $k$-planar graphs that are either $C_k$ free for $k>4$ or are of girth $r$ with $r>5$.
\item Even though we focused on $k$-planar graphs, we believe that extending the study to other beyond-planar graph classes is a challenging research direction that worth to follow. 
\item On the algorithmic side, the recognition problem is of interest; in particular, assuming optimality. A concrete question here, e.g., is whether the problem of recognizing if a graph is optimal $C_3$-free $1$-planar can be done in polynomial time. Recall that in the general setting this problem can be solved in linear time~\cite{DBLP:journals/algorithmica/Brandenburg18}. 
\end{itemize}

\bibliographystyle{plainurl}
\bibliography{k_planar_without_short_cycles_arxiv}
\end{document}